\documentclass[12pt,a4paper,leqno]{amsart}

\usepackage[a4paper,lmargin=2cm,rmargin=2cm,tmargin=4cm,bmargin=4cm]{geometry}

\usepackage[english]{babel}

\usepackage[centertags]{amsmath}
\usepackage{amsfonts}
\usepackage{amssymb}
\usepackage{amsthm}
\usepackage{dsfont}

\usepackage[dvipsnames]{xcolor}

\usepackage[normalem]{ulem}
\usepackage{bbm}

\usepackage[breaklinks=true,colorlinks=true,
linkcolor=MidnightBlue,citecolor=MidnightBlue,
urlcolor=MidnightBlue]{hyperref}

\usepackage{enumerate}

\theoremstyle{plain}

\newtheorem{theorem}{Theorem}

\newtheorem{lemma}{Lemma}

\newtheorem{observation}[lemma]{Observation}

\newtheorem{proposition}[theorem]{Proposition}
\newtheorem{openproblem}[theorem]{Open Problem}

\theoremstyle{remark}
\newtheorem{remark}[theorem]{Remark}
\newtheorem{step}{Step}

\newcommand{\N}{\ensuremath{\mathbb{N}}}

\newcommand{\R}{\ensuremath{\mathbb{R}}}

\newcommand{\Q}{\ensuremath{\mathbb{Q}}}
\newcommand{\eps}{\ensuremath{\varepsilon}}

\newcommand{\NA}{\operatorname{NA}}

\DeclareMathOperator{\sign}{sign}
\DeclareMathOperator{\dist}{dist}
\DeclareMathOperator{\lin}{span}
 
 \newcommand{\tri}[1]{{\left\vert\kern-0.25ex\left\vert\kern-0.25ex\left\vert #1 
    \right\vert\kern-0.25ex\right\vert\kern-0.25ex\right\vert}}

\renewcommand{\leq}{\leqslant}
\renewcommand{\geq}{\geqslant}

\newcommand{\X}{\mathfrak{X}}

\begin{document}
\title[Space with algebraically trivial set of norm-attaining functionals]{A Banach space whose set of norm-attaining functionals is algebraically trivial}

\author[M.\ Mart\'{\i}n]{Miguel Mart\'{\i}n}
\address{Department of Mathematical Analysis and Institute of Mathematics (IMAG), University of Granada, E-18071 Granada, Spain \\
\href{http://orcid.org/0000-0003-4502-798X}{ORCID: \texttt{0000-0003-4502-798X}}}
\email{mmartins@ugr.es}
\urladdr{\url{https://www.ugr.es/local/mmartins}}

\date{July 23rd, 2024; Revised: December 12th, 2024}

\dedicatory{To Mariola, Miguel, and Natalia}

\begin{abstract}
We construct a Banach space $\X$ for which the set of norm-attaining functionals $\NA(\X,\R)$ does not contain any non-trivial cone. Even more, given two linearly independent norm-attaining functionals on $\X$, no other element of the segment between them attains its norm. Equivalently, the intersection of $\NA(\X,\R)$ with a two-dimensional subspace of $\X^*$ is contained in the union of two lines. In terms of proximinality, we show that for every closed subspace $M$ of $\X$ of codimension two, at most four elements of the unit sphere of $\X/M$ have a representative of norm-one. We further relate this example with an open problem on norm-attaining operators.
\end{abstract}


\maketitle

\section{Introduction}
For a real Banach space $X$, we write $X^*$ for its topological dual space (endowed with the usual dual norm), and $B_X$, $S_X$ for its closed unit ball and its unit sphere, respectively. A (closed) subspace $M$ of $X$ is said to be \emph{proximinal} if for every $x\in X\setminus M$, the set
\begin{equation}\label{eq-PM(x)}
P_M(x):=\{m\in M\colon \|x-m\|=\dist(x,M)\}
\end{equation}
is non-empty. This is equivalent to the fact that the quotient map $\pi_M\colon X\longrightarrow X/M$ satisfies $S_{X/M}\subseteq \pi_M(S_X)$. Clearly, every finite-dimensional subspace of $X$ is proximinal, and the same is true for subspaces which are reflexive. On the other hand, a hyperplane is proximinal if and only if it is the kernel of a norm-attaining functional. Recall that a functional $f\in X^*$ is \emph{norm-attaining} if there exists $x\in S_X$ such that $\|f\|=|f(x)|$; we write $\NA(X,\R)$ to denote the set of norm-attaining functionals on $X$. Observe that the Hahn-Banach theorem assures the existence of proximinal hyperplanes in every Banach space and, actually, the Bishop--Phelps theorem assures their abundance.

In the early 1970's, Ivan Singer asked whether every (infinite dimensional) Banach space contains a proximinal subspace of codimension two (see \cite[Problem~1]{Singer-Roumanian}, \cite[Problem~2.1 in p.~14]{Singer}). This question was answered negatively by the late Charles J.\ Read \cite{Read}, who constructed an ingenious renorming $\mathcal{R}$ of $c_0$ without proximinal subspaces of codimension two, and so of any greater finite codimension. Singer's problem is related to the following question by Gilles Godefroy of 2001 \cite[Problem~III]{Godefroy} or \cite[Question~2.26]{Band-Godefroy}: is it true that for every infinite-dimensional Banach space $X$,  $\NA(X,\R)$ contains two-dimensional subspaces? It is immediate that if $Y$ is a proximinal subspace of $X$ of finite-codimension, then $Y^\perp\subset \NA(X,\R)$ (see \cite[Lemma~2.2]{Band-Godefroy}, for instance), but the reversed result is not true in general (see \cite[Section 2]{InduPLMS1982} for instance). In \cite[Theorem~4.2]{Rmoutil}, Martin Rmoutil demonstrates that Read's space $\mathcal R$ also gives a negative solution to Godefroy's question. We refer the reader to \cite{KadetsLopezMartin} to find geometric properties of the space $\mathcal{R}$, to \cite{KLMW-Jussieu} to find further examples of spaces for which the set of norm-attaining functionals contains no two-dimensional subspace, and to \cite{Martin-nbyn} to find $n$-by-$n$ versions of these results. We do not know of other examples constructed in the literature which answer in the negative Singer's question or Godefroy's question. 

Let us comment that both the space $\mathcal{R}$ and the main family of spaces constructed in \cite{KLMW-Jussieu} are not smooth. This implies that if $X$ is any of these spaces, there is an element $x\in S_X$ and two linearly independent elements $f_1,f_2\in S_{X^*}$ such that $f_1(x)=f_2(x)=1$. Hence, $\NA(X,\R)\cap S_{X^*}$ contains a non-trivial segment and, therefore, $\NA(X,\R)$ contains non-trivial cones. There is a construction given in the proof of \cite[Theorem~9.(4)]{DebsGodefroySaint} (see also \cite[Lemma~11]{KadetsLopezMartin} and \cite[Proposition~2.3]{GMZ}) which provides an equivalent smooth norm $|\cdot|$ on a given separable Banach space $X$ such that $\NA((X,|\cdot|),\R)= \NA(X,\R)$. This allows to get ``Read spaces'' which are smooth \emph{a posteriori}. Then, $\NA((X,|\cdot|),\R)\cap S_{(X,|\cdot|),\R)^*}$ does not contain non-trivial segments. However, $\NA((X,|\cdot|),\R)$ contains non-trivial cones. 

There is a somehow related open problem from the 1960's of whether finite-rank operators can be always approximated by norm-attaining ones (see \cite[Question~5]{JohnsonWolfe}, for instance). It is not even known whether there is a Banach space $X$ for which every non-zero norm-attaining operator from $X$ into the two dimensional Hilbert space $\ell_2^2$ is of rank-one. Here is the connection with the above discussion: it is shown in \cite[Theorem~13.9]{KLMW-Lomonosov} that the existence of non-trivial cones in $\NA(X,\R)$ allows to construct a surjective norm-attaining operator from $X$ into $\ell_2^2$. It is natural to ask, as it was done in \cite[Problem~13.11]{KLMW-Lomonosov}, whether $\NA(X,\R)$ contains non-trivial cones for every Banach space $X$.

The aim of this paper is to construct a Banach space $\X$ (which is actually a smooth renorming of $c_0$) for which $\NA(\X,\R)$ contains no non-trivial cones. Moreover, the intersection of $\NA(\X,\R)$ with any two-dimensional subspace $Z$ of $\X^*$ is contained in the union of two one-dimensional subspaces (that is, the intersection of $S_Z$ with $\NA(\X,\R)$ contains at most four points); equivalently, given two linearly independent elements in $\NA(\X,\R)$, no element of the interior of the segment between them belongs to $\NA(\X,\R)$ (all of these justify the name of algebraically trivial in the title of the paper). This is our Theorem~\ref{theorem-main}. We also provide some remarks and comments on it: 
\begin{enumerate}[-]
  \item The norm of $\X$ can be made $\eps$-equivalent to the given norm of $c_0$ for every $\eps>0$ (Remark~\ref{remark-eps-equivalent}).
  \item $\X^{**}$ is strictly convex, hence $\X^*$ is smooth and $\X$ is also strictly convex (Proposition~\ref{Prop-bidualsc}).
  \item The construction in $c_0$ can be carried out in every weakly compactly generated (in particular, separable) Banach space which contains an isomorphic copy of $c_0$ and is (isomorphic to) a closed subspace of $\ell_\infty$ (Proposition~\ref{prop-WCG}). 
\end{enumerate}

We next deduce some proximinality results of the constructed space $\X$ (see Proposition~\ref{prop-proximinal}): given a closed subspace $M$ of $\X$ of codimension two, the set $\{x\in B_\X\colon \pi_M(x)\in S_{\X/M}\}$ contains at most two elements and their opposites, hence $\pi_M(S_\X)\cap S_{\X/M}$ contains at most two elements of $\X/M$ and their opposites. Also, the set 
$$
\{\pi_M(x)\colon x\in \X\setminus M,\ P_M(x)\text{ is not empty}\}
$$
is contained in the union of two one-dimensional subspaces of $\X/M$.

Finally, Subsection~\ref{subsec:NA-operators} contains an explanation of the relation of the constructed space $\X$ and the problem of the denseness of norm-attaining finite-rank operators which was briefly explained above.

\section{The example} \label{section:proofs}
This is the main result of the paper.

\begin{theorem}\label{theorem-main}
There is a Banach space $\X$ such that for every two-dimensional subspace $Z$ of $\X^*$, the intersection of $Z$ with $\NA(\X,\R)$ is contained in the union of two one-dimensional subspaces (that is, the intersection of $S_Z$ with $\NA(\X,\R)$ contains at most four points). Equivalently, given two linearly independent elements in $\NA(\X,\R)$, no element of the interior of the segment between them belongs to $\NA(\X,\R)$.
\end{theorem}

To prove the theorem, we will need a series of preliminary results, some of which may be of independent interest. We start constructing equivalent norms on $\ell_1$ which are smooth but still share some properties of the usual duality between $\ell_1$ and $\ell_\infty$.

\begin{lemma}\label{lemma-A}
Let $\Phi$ be a sequence of positive real numbers which belongs to $\ell_1$ with $\|\Phi\|_1<1$, and consider the operator $S_\Phi\colon \ell_2 \longrightarrow \ell_1$ by $[S_\Phi(a)](n)=\Phi(n)a(n)$ for every $n\in \N$ and every $a\in \ell_2$. We define an equivalent norm $|\cdot|_\Phi$ on $\ell_1$ such that
$$
B_{(\ell_1,|\cdot|_\Phi)}=B_{\ell_1} + S_\Phi(B_{\ell_2}).
$$
Then:
\begin{enumerate}[(a)]
\item for $f\in \ell_\infty$, $|f|_\Phi^*=\|f\|_\infty + \|S_\Phi^*(f)\|_2$;
\item for $f\in \ell_\infty$, $|f|_\Phi^*\geq \|f\|_\infty \geq (1 -\|\Phi\|_1)|f|_\Phi^*$;
\item for $x\in \ell_1$, $|x|_\Phi\leq \|x\|_1\leq (1+\|\Phi\|_1)|x|_\Phi$;
\item $(\ell_1,|\cdot|_\Phi)^*$ is strictly convex, so $(\ell_1,|\cdot|_\Phi)$ is smooth;
\item if $x\in \ell_1$ and $f\in \ell_\infty$ satisfy that $|f|_\Phi^*=1$, $\langle f,x\rangle = |x|_\Phi$, and $|x(n)|>\Phi(n)|x|_\Phi$, it follows that $f(n)=\sign (x(n))\|f\|_\infty$.
\end{enumerate}
\end{lemma}

\begin{proof}
We follow the arguments in \cite[Theorem 9.(4)]{DebsGodefroySaint} (given also at \cite[Lemma~11]{KadetsLopezMartin}) and the more abstract version which appears in \cite[Lemma~2.2]{KLMW-Jussieu}. First, the existence of the norm $|\cdot|_\Phi$ and the formula for the dual norm in (a) follows directly from \cite[Lemma~2.2]{KLMW-Jussieu}. 
To get (b), just observe that $\|S_\Phi\|\leq \|\Phi\|_1$ and then 
$$
\|S_\Phi^*(f)\|_2\leq \|S_\Phi\|\|f\|_\infty \leq \|S_\Phi\|\,|f|_\Phi^*.
$$
(c). $B_{\ell_1}\subset B_{(\ell_1,|\cdot|_\Phi)}\subset (1+\|S_\Phi\|)B_{\ell_1}\subset (1+\|\Phi\|_1)B_{\ell_1}$.

(d) follows from the fact that $S_\Phi^*$ is one-to-one and the norm of $\ell_2$ is strictly convex using the formula given in (a).

(e). Suppose that $f\in S_{(\ell_1,|\cdot|_\Phi)^*}$ attains its norm at $x\in S_{(\ell_1,|\cdot|_\Phi)}$. Then, there is $x_0\in S_{\ell_1}$ and $a_0\in S_{\ell_2}$ such that $f$ attains its maximum on $B_{\ell_1}$ (equals to $\|f\|_\infty$) at $x_0$ and its maximum on $S_\Phi(B_{\ell_2})$ (equals to $\|S_\Phi^*(f)\|_2$) at $S_\Phi(a_0)$, and $x=x_0+S_\Phi(a_0)$. Now, for those $n\in \N$ such that $x_0(n)\neq 0$, it follows that $f(n)=\sign(x_0(n))\|f\|_\infty$. But if $x(n)>\Phi(n)$, then $x_0(n)\geq x(n)- \Phi(n)|a_0(n)|>\Phi(n)(1-|a_0(n)|)\geq 0$, hence $f(n)=\|f\|_\infty$. Analogously, if $x(n)<-\Phi(n)$, we get $f(n)=-\|f\|_\infty$. Finally, the result for $x$ with arbitrary norm follows by homogeneity.
\end{proof}

Next, we adapt Proposition~2.8 of \cite{KLMW-Jussieu} to get sequences which are ``doubly infinitely dense'' through an injective operator in any operator range. Recall that a linear subspace $Y$ of a Banach space $X$ is said to be an \emph{operator range} if there is an infinite-dimensional Banach space $E$ and a bounded injective operator $T\colon  E \longrightarrow X$ such that $T(E) = Y$. Equivalently, there is a complete norm on $Y$ which is stronger than the given norm of $Y$. We refer to the papers \cite{Cross,CrossOstrovskiiShevchik, JimenezLajara} and references therein for more information and background.

\begin{lemma}\label{lemma-B}
Let $Y$ be a separable operator range. Then, there exists a norm-one injective operator $T\colon \ell_1(\N\times \N) \longrightarrow Y$ such that for every $m\in \N$, the set
$$
\left\{\frac{T(e_{n,m})}{\|T(e_{n,m})\|}\colon n\in \N\right\}
$$  
is dense in $S_Y$.
\end{lemma}

\begin{proof}
Consider $\sigma\colon\N\times \N\longrightarrow \N$ bijective and the isometric isomorphism $\Psi\colon\ell_1(\N\times\N) \longrightarrow \ell_1$ which carries $e_{n,m}\in \ell_1(\N\times\N)$ to $e_{\sigma(n,m)}\in \ell_1$ for every $n,m\in \N$.

Now, take a sequence $\{w'_m\colon m\in \N\}$ such that every element of $S_Y$ is an accumulation point of it and repeat the proof of Proposition 2.8 of \cite{KLMW-Jussieu} with the sequence $w_{\sigma(n,m)}=w'_n$ for every $n,m\in \N$.
\end{proof}

We will also make use of the following two easy observations.

\begin{observation}\label{lemma-C}
Let $\{Z_m\colon m\in \N\}$ be a family of smooth Banach spaces and let $Z=\left[\bigoplus_{m\in \N} Z_m\right]_{\ell_1}$ be its $\ell_1$ sum. If $z=(z_m)\in Z$ and $f=(f_m)\in Z^*= \left[\bigoplus_{m\in \N} Z_m^*\right]_{\ell_\infty}$ satisfy that $\|f\|_\infty=1$ and
$\langle f,z\rangle =\|z\|_1$, then $f_m(z_m)=\|z_m\|$ for every $m\in \N$, hence $\|f_m\|=1$ when $z_m\neq 0$. As a consequence, the norm of $Z$ is smooth at every element $z=(z_m)\in X$ such that $z_m\neq 0$ for every $m\in \N$.  
\end{observation}

\begin{observation}\label{observation-D}
Let $Z$ be a smooth Banach space and let $z_1^*,z_2^*\in S_{Z^*}$ be linearly independent. If for some $0<t<1$ we have $\|tz_1^*+(1-t)z_2^*\|=1$, then $tz_1^*+(1-t)z_2^* \notin \NA(Z,\R)$.
\end{observation}

We are now ready to prove the main result of the paper.

\begin{proof}[Proof of Theorem~\ref{theorem-main}]
We divide the job in ten steps for clarity of the exposition.

\begin{step}\label{step-1}
Consider $X=c_0$ endowed with an equivalent norm $\tri{\cdot}$ which is smooth and satisfies that $\NA((X,\tri{\cdot}),\R)=\NA(c_0,\R)$, see \cite[Lemma~11]{KadetsLopezMartin} or the proof of \cite[Theorem~9.(4)]{DebsGodefroySaint}. It then follows from \cite[Proposition~2.4]{KLMW-Jussieu} that there is a separable operator range $Y$ dense in $(X,\tri{\cdot})^*$ with 
$$
\NA((X,\tri{\cdot}),\R)\cap Y=\{0\}.
$$
\end{step}

\begin{step}\label{step-2}
By Lemma~\ref{lemma-B}, there exists a norm-one operator $T\colon \ell_1(\N\times \N)\longrightarrow Y\subset (X^*,\tri{\cdot})$ which is injective and satisfies that for every $m\in \N$, the set
\begin{equation}\label{eq-step2-set}
\left\{\frac{T(e_{n,m})}{\tri{T(e_{n,m})}}\colon n\in \N\right\}
\end{equation}
is dense in $S_Y$ (hence, dense in $S_{(X,\tri{\cdot})^*}$). Write $v^*_{n,m}:=T(e_{n,m})\in B_Y\subset B_{(X,\tri{\cdot})^*}$ for every $n,m\in \N$.
\end{step}

\begin{step}\label{step-3}
Fix $m\in \N$. Let $|\cdot|_m$ be the equivalent norm on $\ell_1$ given by Lemma~\ref{lemma-A} for the sequence $\displaystyle \Phi_m(n)=\frac{1}{2^m}\frac{1}{2^n}\tri{v_{n,m}^*}$ for every $n\in \N$. Then:
\begin{enumerate}[(i)]
\item $(\ell_1,|\cdot|_m)^*$ is strictly convex, so $(\ell_1,|\cdot|_m)$ is smooth;
\item for $x\in \ell_1$, $|x|_m\leq \|x\|_1\leq \left(1+\frac{1}{2^m}\right)|x|_m$; 
\item for $f\in (\ell_1,|\cdot|_m)^*$, $|f(n)|\leq |f|_m^*$ for every $n\in \N$;
\item if $x\in (\ell_1,|\cdot|_m)$, $f\in (\ell_1,|\cdot|_m)^*$, and $n\in \N$ satisfy that 
    $$
    |f|^*_m=1,\quad \langle f,x\rangle = |x|_m,\quad |x(n)|> \frac{1}{2^m}\frac{1}{2^n}\tri{v_{n,m}^*}|x|_m,
    $$
then
    $$
    f(n)=\sign(x(n))\|f\|_\infty \quad \text{and} \quad \|f\|_\infty\in \left[1-\frac{1}{2^m},1\right].
    $$
\end{enumerate}
\end{step}

\begin{step}\label{step-3ymedio}
Define $R_m\colon (X,\tri{\cdot})\longrightarrow (\ell_1,|\cdot|_m)$ by 
$$
[R_m (x)](n)=\frac{m}{2^m}\frac{1}{2^n}v_{n,m}^*(x) \qquad \bigl(n\in \N,\ x\in X\bigr).
$$
Then:
\begin{enumerate}[(a)]
  \item $R_m$ is one-to-one since $\{v_{n,m}^*\colon n\in \N\}$ separates the points of $(X,\tri{\cdot})$.
  \item $|R_m(x)|_m\leq \|R_m(x)\|_1\leq \frac{m}{2^m}\tri{x}$ for every $x\in X$.
\end{enumerate}
\end{step}

\begin{step}\label{step-4}
Define $R\colon (X,\tri{\cdot})\longrightarrow W:=\left[\bigoplus\limits_{m\in \N} (\ell_1,|\cdot|_m)\right]_{\ell_1}$ by 
$$
R(x) = \bigl(R_m (x)\bigr)_{m\in \N} \qquad (x\in X).
$$
Then, for $x\in X$, 
$$
\|R(x)\|_W= \sum_{m=1}^{\infty} |R_m (x)|_m \leq \left(\sum_{m=1}^{\infty}\frac{m}{2^m}\right)\tri{x}.
$$
By Step~\ref{step-3ymedio}.(a), given $x\in X\setminus \{0\}$, all coordinates of $R(x)=\bigl(R_m (x)\bigr)_{m\in \N}$ are different from zero. This shows that $R$ is one-to-one and also that the norm of $W$ is smooth at $R(x)$ for every $x\in X\setminus \{0\}$ by Observation~\ref{lemma-C}.
\end{step}

\begin{step}\label{step-5}
Given $(w_m^*)_{m\in \N}\in W^*\equiv \left[\bigoplus\limits_{m\in \N} (\ell_1,|\cdot|_m)^*\right]_{\ell_\infty}$, one has that
\begin{equation}\label{eq:step4--eq1-R*inY}
R^*\bigl((w_m^*)_{m\in\N} \bigr)=T\left(\sum_{m=1}^{\infty}\sum_{n=1}^{\infty} w_m^*(n)\frac{m}{2^m}\frac{1}{2^n} e_{n,m} \right) \in Y\subset (X,\tri{\cdot})^*,
\end{equation}
where $T\colon \ell_1(\N\times \N)\longrightarrow Y$ is the operator considered in Step~\ref{step-2}; moreover, 
\begin{equation}\label{eq:step-4--eq2-allzero}
R^*\bigl((w_m^*)_{m\in\N} \bigr)=0 \quad \Longrightarrow \quad w_m^*(n)=0 \ \text{ for all $n,m\in \N$}.
\end{equation}
Indeed, Equation~\eqref{eq:step4--eq1-R*inY} follows routinely from the fact that the series
$$
\sum_{m=1}^{\infty}\sum_{n=1}^{\infty} w_m^*(n)\frac{m}{2^m}\frac{1}{2^n} e_{n,m}
$$
is absolutely summable in $\ell_1$ as $|w_m^*(n)|\leq |w_m^*|_m^*\leq \left\|(w_m^*)_{m\in \N}\right\|_\infty$ (by (iii) in Step~\ref{step-3}).
As $T$ is one-to-one and $\{e_{n,m}\colon n\in \N,\, m\in \N\}$ is a basis of $\ell_1(\N\times \N)$, Equation~\eqref{eq:step-4--eq2-allzero} follows immediately.
\end{step}

\begin{step}\label{step-6}
Define an equivalent norm $p$ on $X$ by
$$
p(x)=\tri{x} + \|R(x)\|_W \qquad (x\in X)
$$
and consider $\X=(X,p)$. Then, $\X$ is smooth since $\tri{\cdot}$ is a smooth norm and the map $x\longmapsto \|Rx\|_W$ is smooth at every non-zero element of $X$, see Step~\ref{step-4}. 

By Lemma~2.1.c in \cite{KLMW-Jussieu}, if $x\in \X$ and $x^*\in \X^*$ satisfy that $p(x)=p(x^*)=x^*(x)=1$, then
there is $x^*_0\in \NA(X,\tri{\cdot})$ with $\tri{x_0^*}=1$ and $w^*:=(w_m^*)_{m\in \N}\in S_{W^*}$ such that
\begin{equation*}
x^*=x_0^* + R^*(w^*) \quad \text{and} \qquad \langle (w_m^*)_{m\in \N},R(x)\rangle=\|R(x)\|_W.
\end{equation*}
Moreover, since $R_m(x)\neq 0$ for every $m\in \N$ (see Step~\ref{step-3ymedio}.(a)),  Observation~\ref{lemma-C} gives that
\begin{equation*}
|w_m^*|^*_m=1 \quad  \text{and} \quad w_m^*(R_m x)=|R_m x|_m \qquad (m\in \N).
\end{equation*}
\end{step}

\begin{step}\label{step-7}
Pick $x^*,z^*\in \NA(\X,\R)$ with $p(x^*)= p(z^*)= 1$ which are linearly independent and that attain their norms at $x,z\in S_\X$, respectively. Hence, $x,z$ are linearly independent as the norm of $\X$ is smooth. Suppose, for the sake of getting a contradiction, that $u^*:=tx^*+(1-t)z^*\in \NA(\X,\R)$ for some $0<t<1$ and that it attains its norm at $u\in S_\X$. Then, by Step~\ref{step-6}, there are $x_0^*,z_0^*,u_0^*$ in $S_{(X,\tri{\cdot})^*}\cap \NA(X,\tri{\cdot})$, sequences $\Xi=(\xi_m^*)_{m\in \N}$, $\Theta=(\zeta_m^*)_{m\in \N}$, $\Upsilon=(u_m^*)_{m\in \N}$ in $S_{W^*}$ satisfying
\begin{equation*}
x^*=x_0^* + R^*(\Xi),\quad z^*=z_0^* + R^*(\Theta),\quad \frac{tx^*+(1-t)z^*}{p(tx^*+(1-t)z^*)}=u_0^* + R^*(\Upsilon),
\end{equation*}
and
\begin{equation}\label{eq:step-7-eq2}
\xi_m^*(R_m x)=|R_m x|_m, \quad \zeta_m^*(R_m z)=|R_m z|_m, \quad u_m^*(R_m u)=|R_m u|_m \qquad (m\in \N).
\end{equation}
Therefore, 
$$
t x_0^* + (1-t)z_0^* - p(tx^*+(1-t)z^*)u_0^* = - R^*\bigl(t\,\Xi + (1-t) \Theta - p(tx^*+(1-t)z^*)\Upsilon\bigr).
$$
As the left-hand side belongs to $\lin \NA(X,\tri{\cdot})=\NA(X,\tri{\cdot})$, the right-hand side belongs to $Y$ by Step~\ref{step-5}, and $\NA(X,\tri{\cdot})\cap Y=\{0\}$ by Step~\ref{step-1}, it follows that 
$$
R^*\bigl(t\,\Xi + (1-t) \Theta - p(tx^*+(1-t)z^*)\Upsilon\bigr)=0.
$$ 
Using now Equation~\eqref{eq:step-4--eq2-allzero}, we get that 
\begin{equation}\label{eq-step7-eq-main}
t\, \xi_m^*(n) + (1-t) \zeta^*_m(n) = p(tx^*+(1-t)z^*)u_m^*(n) \quad \text{for all $n,m\in \N$}.
\end{equation}
\end{step}

\begin{step}\label{step-8}
As $x,z$ are linearly independent, we may find $\phi\in S_{(X,\tri{\cdot})^*}$ and $\delta>0$ such that 
$$
\phi(x)>\delta,\qquad \phi(z)>\delta.
$$
Indeed, just observe that $x$ does not belong to the one-dimensional subspace $\R(x-z)$, and Hahn-Banach theorem gives $\phi\in S_{(X,\tri{\cdot})^*}$ such that $\phi(x)=\dist(x,\R(x-z))$ and $\phi(x-z)=0$.
\end{step}

\begin{step}\label{step-9}
Consider $m\in \N$ such that 
\begin{equation}\label{eq-step9-eq1}
\frac{1}{m} < \delta \quad \text{and} \quad p(tx^*+(1-t)z^*)<1-\frac{1}{2^m}
\end{equation}
(the second inequality is possible as Observation~\ref{observation-D} gives that $p(tx^*+(1-t)z^*)<1$ since $tx^*+(1-t)z^*\in \NA(\X,\R)$). Now, as the set
$$\left\{\frac{v^*_{n,m}}{\tri{v^*_{n,m}}}\colon n\in \N\right\}$$ is dense in $S_Y$ by Step~\ref{step-2}, hence dense in $S_{(X,\tri{\cdot})^*}$, we may find $n\in \N$ such that 
\begin{equation*}
v_{n,m}^*(x)>\delta \tri{v_{n,m}^*}  \quad \text{and} \quad v_{n,m}^*(z)>\delta \tri{v_{n,m}^*}.
\end{equation*}
Then, 
\begin{align*}
  [R_m x](n) & = \frac{m}{2^m}\frac{1}{2^n}v_{n,m}^*(x)> \frac{m}{2^m}\frac{1}{2^n}\delta \tri{v_{n,m}^*}  > \frac{1}{2^m}\frac{1}{2^n} \tri{v_{n,m}^*}\geq \frac{1}{2^m}\frac{1}{2^n} \tri{v_{n,m}^*}|R_m x|_m.
\end{align*}
This together with Equation~\eqref{eq:step-7-eq2} imply by Step~\ref{step-3}.(iv) that
$$
\xi_m^*(n)=\|\xi_m^*\|_\infty \in \left[1-\frac{1}{2^m},1\right].
$$
Analogously, 
$$
\zeta_m^*(n)=\|\zeta_m^*\|_\infty \in \left[1-\frac{1}{2^m},1\right].
$$
By Equation~\eqref{eq-step7-eq-main}, we have for the fixed $n,m\in \N$ that
$$
t\, \xi_m^*(n) + (1-t) \zeta^*_m(n) = p(tx^*+(1-t)z^*)u_m^*(n).
$$
But the left-hand side belongs to the interval $\left[1-\frac{1}{2^m},1\right]$, while the right-hand side has absolute value strictly smaller than $1-\frac{1}{2^m}$ by Equation~\eqref{eq-step9-eq1}, a contradiction.\qedhere
\end{step}
\end{proof}

\subsection{Remarks and further results on the construction}
Let us gives some comments on the construction made above.

First, the space $\X=(X,p)$ constructed in the proof of Theorem~\ref{theorem-main} is isomorphic to $c_0$ and there is no problem in getting $p$ to be $\eps$-equivalent to the given norm of $c_0$. We say that a norm $p$ on a Banach space $X$ is \emph{$\eps$-equivalent} to the given norm $\|\cdot\|$ if there are $M,m>0$ with $M/m\leq 1+\eps$, such that 
$$
m\|x\|\leq p(x)\leq M\|x\| \qquad (x\in X).
$$

\begin{remark}\label{remark-eps-equivalent}
{\slshape Given $\eps>0$, there is a norm $p_\eps$ on $c_0$, which is $\eps$-equivalent to the given norm of $c_0$, such that $\X_\eps=(c_0,p_\eps)$ satisfies the properties of the space $\X$ in Theorem~\ref{theorem-main}.} Indeed, for $\delta>0$, consider the norm $\tri{\cdot}$ in $c_0$ of Step~\ref{step-1} to be $\delta$-equivalent to the given norm of $c_0$ (this can be easily done using \cite[Lemma~11]{KadetsLopezMartin});  in Step~\ref{step-6}, consider $p_\delta(x)=\tri{x} + \delta\|Rx\|_W$ for every $x\in X$. It readily follows that $p_\delta$ is $\eps$-equivalent to $\|\cdot\|_\infty$ of $c_0$ for sufficiently small $\delta$.
\end{remark}

We next show some geometric properties of the space $\X$.

\begin{proposition}\label{Prop-bidualsc}
The space $\X$ constructed in the proof of Theorem~\ref{theorem-main} satisfies that $\X^{**}$ is strictly convex, hence $\X^*$ is smooth.
\end{proposition}

\begin{proof}
As all the operators $R_m$ with $m\in \N$ and $R$ given in Steps \ref{step-3ymedio} and \ref{step-4} are compact, we have
$$
R^{**}(\X^{**})\subseteq 
J_W(W)\equiv W \quad\text{and}\quad R_m^{**}(\X^{**})\subseteq 
J_{(\ell_1,|\cdot|_m)}(\ell_1,|\cdot|_m)\equiv (\ell_1,|\cdot|_m) \ \ (m\in \N).
$$
Hence, the bidual norm of the norm $p$ of $\X$ given in Step~\ref{step-6} has the formula 
\begin{equation}\label{eq-bidual-sc-1}
p(x^{**})=\tri{x^{**}} + \|R^{**}(x^{**})\|_{W} = \tri{x^{**}} + \sum_{m=1}^{\infty} |R_m^{**}(x^{**})|_m 
 \qquad (x^{**}\in \X^{**})
\end{equation}
by \cite[Lemma~2.1(b)]{KLMW-Jussieu}. Moreover, the formula given in Step~\ref{step-3ymedio} for $R_m$ extends to $R_m^{**}$, that is, for $x^{**}\in \X^{**}$, $R_m^{**}(x^{**})\in (\ell_1,|\cdot|_m)$ with coordinates
$$
[R^{**}_m (x^{**})](n)=\frac{m}{2^m}\frac{1}{2^n}x^{**}(v_{n,m}^*) \qquad \bigl(n\in \N, \ m\in \N).
$$ 
Suppose, for the sake of contradiction, that two linearly independent elements $x^{**}, z^{**}\in \X^{**}$ satisfy that
$$
p(x^{**}+z^{**})=p(x^{**})+p(z^{**}).
$$
It then follows from Equation~\eqref{eq-bidual-sc-1} that
\begin{equation*}
\bigl|R_m^{**}(x^{**}+z^{**})\bigr|_m = |R_m^{**}(x^{**})|_m + |R_m^{**}(z^{**})|_m \ \ \ (m\in\N).
\end{equation*}
As $x^{**}, z^{**}\in \X^{**}$ are linearly independent, by the norm density of $S_Y$ in $S_{(X,\tri{\cdot})^*}$, we may find $y^*\in S_Y$ and $m\in \N$ such that
$$
x^{**}(y^*)<-\frac{1}{m}<0<\frac{1}{m}<z^{**}(y^*).
$$
Pick $f\in (\ell_1,|\cdot|_m)^*$ with $|f|_m^*=1$ satisfying that 
$$
\langle f,R_m^{**}(x^{**}+z^{**})\rangle=\bigl|R_m^{**}(x^{**}+z^{**})\bigr|_m = |R_m^{**}(x^{**})|_m + |R_m^{**}(z^{**})|_m
$$
(Recall that $R_m^{**}(\X^{**})\subseteq (\ell_1,|\cdot|_m)$.) 
It then follows that
\begin{equation}\label{eq:bidualsc}
\langle f, R_m^{**}(x^{**})\rangle= |R_m^{**}(x^{**})|_m \qquad \text{and} \qquad \langle f,R_m^{**}(z^{**})\rangle=|R_m^{**}(z^{**})|_m.
\end{equation}
On the other hand, as the set
$
\left\{\dfrac{v^*_{n,m}}{\tri{v^*_{n,m}}}\colon n\in \N\right\}
$
is norm dense in $S_Y$ by Step~\ref{step-2}, we may find $n\in \N$ such that 
\begin{equation*}
x^{**}(v_{n,m}^*)<-\frac{1}{m}\tri{v_{n,m}^*}  \qquad \text{and} \qquad z^{**}(v_{n,m}^*)>\frac{1}{m}\tri{v_{n,m}^*}.
\end{equation*}
Then, 
$$
  [R_m^{**} x^{**}](n)  = \frac{m}{2^m}\frac{1}{2^n}x^{**}(v^*_{n,m})< \frac{-1}{2^m}\frac{1}{2^n} \tri{v_{n,m}^*}\leq \frac{-1}{2^m}\frac{1}{2^n} \tri{v_{n,m}^*}\,|R_m^{**} x^{**}|_m
$$
and
$$
   [R_m^{**} z^{**}](n)  = \frac{m}{2^m}\frac{1}{2^n}z^{**}(v_{n,m}^*)> \frac{1}{2^m}\frac{1}{2^n} \tri{v_{n,m}^*}\geq \frac{1}{2^m}\frac{1}{2^n} \tri{v_{n,m}^*}\,|R_m^{**} z^{**}|_m.
$$
This, together with Equation~\eqref{eq:bidualsc}, imply by Step~\ref{step-3}.(iv) that
$$
f(n)=-\|f\|_\infty<0 \qquad \text{and} \qquad f(n)=\|f\|_\infty>0,
$$
a clear contradiction.
\end{proof}

In a earlier version of this manuscript, there was an error in the proof of the above proposition which was detected by the referee. The same error has also occurred in the
proof of \cite[Proposition~3.4]{KLMW-Jussieu}, but a similar solution of using the norm-density of $\left\{\frac{v_n}{\|v_n\|}\right\}$ instead of their weak-star density works
there as well.

The process done in Theorem~\ref{theorem-main} can be repeated in every weakly compactly generated Banach space $X$ which is isomorphic to a subspace of $\ell_\infty$ and contains an isomorphic copy of $c_0$.

\begin{proposition}\label{prop-WCG}
Let $X$ be a weakly compactly generated (for instance, separable) Banach space which is isomorphic to a subspace of $\ell_\infty$ and contains an isomorphic copy of $c_0$. Then, there is an equivalent norm $p$ on $X$ such that for every two-dimensional subspace $Z$ of $(X,p)^*$, the intersection of $Z$ with $\NA((X,p),\R)$ is contained in the union of two one-dimensional subspaces (that is, the intersection of $S_Z$ with $\NA((X,p),\R)$ contains at most four points). Equivalently, given two linearly independent elements in $\NA((X,p),\R)$, no element of the interior of the segment between them belongs to $\NA((X,p),\R)$.
\end{proposition}

\begin{proof} By \cite[Theorem~4.1]{KLMW-Jussieu}, there is an equivalent norm $\|\cdot\|_1$ on $X$ and a separable operator range $Y$ in $(X^*,\|\cdot\|_1)$ which is weak-star dense in $(X^*,\|\cdot\|)$ and such that $\NA((X,\|\cdot\|_1),\R)\cap Y=\{0\}$. Next, 
\cite[Proposition~2.3]{GMZ} allows us to produce a smooth norm $\tri{\cdot}$ on $X$ such that $\NA((X,\tri{\cdot}),\R)=\NA((X,\|\cdot\|),\R)$. Hence, $(X,\tri{\cdot})$ may play the role of $c_0$ in the proof of Theorem~\ref{theorem-main}, with the only difference that $S_Y$ would not necessarily be dense in $S_{(X,\tri{\cdot})^*}$, but only $Y$ is weak-star dense in $(X,\tri{\cdot})^*$. But this is enough to get Step~\ref{step-3ymedio}.(a) as the set given in Equation~\eqref{eq-step2-set} separates the points of $(X,\tri{\cdot})^*$; in Step~\ref{step-9}, even though $S_Y$ is not necessarily norm dense in $S_{(X,\tri{\cdot})^*}$, the set of restrictions of elements in $S_Y$ to $\lin\{x,z\}$ is dense in $\lin\{x,z\}^*$, the property which is actually used.
\end{proof}

We do not know whether every non-reflexive Banach space can be renormed in such a way that the set of norm-attaining functionals in the new norm fails to contain non-trivial cones. However, having an equivalent norm with the properties of the one of $\X$ in Theorem~\ref{theorem-main} is not possible for Banach spaces with the Radon-Nikod\'{y}m property. Indeed, it is shown in \cite[Proposition~5.6]{KLMW-Jussieu} that, in this case, the set of norm-attaining functionals contains the $\Q$-linear span of an $\R$-linearly independent infinite sequence. In particular, there are two-dimensional subspaces of the dual for which their intersection with the set of norm-attaining functionals is dense in the subspace.

The following remark on a result by Vesel\'{y} of 2009 is pertinent.

\begin{remark}\label{remarkVesely}
{\slshape
It was proved by L.~Vesel\'{y} in \cite[Theorem 3.2]{Vesely} that for every Banach space $Z$, the set $\NA(Z,\R)\cap S_{Z^*}$ is $\boldsymbol{c}$-dense in $S_{Z^*}$ (that is, every open subset of $S_{Z^*}$ contains at least continuum many elements of $\NA(Z,\R)$). Moreover, if $Z$ admits an equivalent Fr\'{e}chet differentiable norm (in particular, if $Z$ is isomorphic to $c_0$ as our space $\X$ in Theorem~\ref{theorem-main}), the set $\NA(Z,\R)\cap S_{Z^*}$ is pathwise connected and locally pathwise connected. These results make the existence of a space like $\X$ even more surprising, if possible.}
\end{remark}

\subsection{Some proximinality properties of the constructed space}\label{subsection-proximinality} It is clear that no two-codimensional subspace of $\X$ is proximinal (as $\NA(\X,\R)$ contains no two-dimensional subspaces). The next result shows that, actually, the situation is even more extreme.

\begin{proposition}\label{prop-proximinal}
Let $M$ be a two-codimensional closed subspace of $\X$. Then:
\begin{enumerate}[(a)]
    \item The set 
    $$
    \bigl\{x\in B_\X\colon \pi_M(x)\in S_{\X/M}\bigr\}
    $$
    contains, at most, two distinct elements and their opposites. 
    \item There are two elements $x_1,x_2\in B_\X$ such that  
      $$
      \pi_M(B_\X)\cap S_{\X/M}\subseteq \{x_1+M, x_2+M, -x_1+M, -x_2+M\}.
      $$
  \item There are two elements $x_1,x_2\in \X$ such that the (maybe empty) set of elements in the quotient $\X/M$ given by
  $$
  \bigl\{x+M\colon x\in \X\setminus M,\ P_M(x) \text{ not empty}\bigr\}
  $$
  is contained in $\R(x_1+M)\cup \R(x_2+M)$.
\end{enumerate}
\end{proposition}

\begin{proof}
(a). Take $x_1,x_2,x_3\in B_\X$ satisfying that $\pi_M(x_i)\in S_{\X/M}$ for $i=1,2,3$. Pick $\phi_1,\phi_2,\phi_3\in (X/M)^*$ with $\phi_i(x_i+M)=\|\phi_i\|=1$ for $i=1,2,3$. Now, for $i=1,2,3$, writing $f_i:=\phi_i\circ \pi_M\in M^\perp\subset \X^*$, we get that $\|f_i\|=f_i(x_i)=1$ so, in particular, $f_i\in \NA(\X, \R)$. As $M^\perp$ is a two-dimensional subspace of $\X^*$, Theorem~\ref{theorem-main} gives that there are $i,j\in \{1,2,3\}$, $i\neq j$, such that $f_i=\pm f_j$. Since $\X$ is strictly convex by  Proposition~\ref{Prop-bidualsc}, we get that $x_i=\pm x_j$. Assertion (b) follows immediately from (a).

(c). Let $x_1,x_2$ given by (b). Let $x\in \X\setminus M$ such that $P_M(x)$ is not empty and take $m\in M$ such that $\|x-m\|=\dist(x,M)=\|x+M\|$. Then, $\dfrac{x-m}{\|x-m\|}\in S_\X$ satisfies that for every $m'\in M$, 
$$
\left\|\dfrac{x-m}{\|x-m\|}-m'\right\|=\frac{\bigl\|x-(m+\|x-m\|m'\bigr\|}{\|x-m\|}\geq \frac{\dist(x,M)}{\|x-m\|}=1,
$$
hence $\pi_M\left(\dfrac{x-m}{\|x-m\|}\right)\in S_{\X/M}$. It follows that 
$$
\pi_M\left(\dfrac{x-m}{\|x-m\|}\right)\in \{\pm (x_1 + M), \pm (x_2 + M)\},
$$
that is, $\pi_M(x)\in \R(x_1+M)\cup \R(x_2+M)$.
\end{proof}

\subsection{Relation with norm-attaining finite-rank operators}\label{subsec:NA-operators}
Given Banach spaces $X$, $Y$, we write $\mathcal{L}(X,Y)$ to denote the Banach space of all bounded linear operators from $X$ to $Y$ endowed with the operator norm. An operator $T\in \mathcal{L}(X,Y)$ is said to be \emph{norm-attaining} (writing $T\in \NA(X,Y)$) if there is $x\in S_X$ such that $\|Tx\|=\|T\|$. Since the seminal paper by Lindenstrauss of 1963 \cite{Lin}, there has been an extensive study of when bounded linear operators can be approximated (in norm) by norm-attaining ones. This is the case if the domain is reflexive (Lindenstrauss), or even if it has the Radon-Nikod\'{y}m property (Bourgain), or when the range is finite-dimensional and polyhedral or it is isometric to a subspace of $\ell_\infty$ containing the canonical copy of $c_0$ (Lindenstrauss). On the other hand, we may find operators which cannot be approximated by norm-attaining ones if either the domain is $c_0$ or $C[0,1]$ or $L_1[0,1]$ (Lindenstrauss) or the range is $\ell_p$ with $1\leq p<\infty$ (Gowers, Acosta). There are also compact operators which cannot be approximated by norm-attaining ones (Mart\'{\i}n). The interested reader may find a lot of information in \cite{Aco-survey,Bachir, JMR-JFA2023,KLMW-Lomonosov, Martin2016}, and references therein.

Let us mention that it is still not known whether finite-rank operators can be always approximated by norm-attaining ones, even in the case when the range space is the two-dimensional Hilbert space $\ell_2^2$ (in the words of Johnson and Wolfe \cite[Question~6]{JohnsonWolfe}, the latter is \emph{the most irritating open problem about norm-attaining operators}).  Actually, it is still unknown if there is a Banach space $X$ such that $\NA(X,\ell_2^2)$ does not contain rank-two operators. Here is where the space $\X$ comes into play in this setting. It is known that a finite rank operator $T\in \mathcal{L}(X,Y)$ belongs to $\NA(X,Y)$ whenever $[\ker T]^\perp\subseteq \NA(X,\R)$ \cite[Corollary~13.6]{KLMW-Lomonosov}. Hence there are elements in $\NA(X,\ell_2^2)$ of rank two whenever $\NA(X,\R)$ contains two-dimensional subspaces. Moreover, it is enough that $\NA(X,\R)$ contains non-trivial cones to get that $\NA(X,Y)$ contains rank-two operators for every Banach space $Y$ of dimension greater than or equal to two \cite[Corollary~13.10]{KLMW-Lomonosov}, and this applies to all previously constructed Read norms, as explained in the introduction, see \cite[Example~13.13]{KLMW-Lomonosov}. However, $\NA(\X,\R)$ does not contain non-trivial cones, hence no known result implies that $\NA(\X,\ell_2^2)$ contains rank-two operators.

\begin{openproblem}
Does $\NA(\X,\ell_2^2)$ contain rank-two operators?
\end{openproblem}

Let us mention that \cite[Theorem~13.24]{KLMW-Lomonosov} contains a very general sufficient condition on $\NA(X,\R)$ (which implies that it contains subspaces of arbitrary large dimension) to assure that every finite-rank operator whose domain is $X$ can be approximated by norm-attaining ones. 

Finally, we would like to draw the interested reader's attention to an open question related to the linear structure of the set of norm-attaining functionals \cite[Question~2.25]{Band-Godefroy}: does there exist a non-reflexive Banach space $X$ for which $\NA(X^*,\R)$ is a vector space? It is known that in such a case, $X$ has to be weakly sequentially complete \cite[Proposition~2.24]{Band-Godefroy}. More related open questions can be found in \cite[\S 3]{GodeInduJAT}. To find applications of norm attainment, the reader is encouraged to consult \cite{Godefroy-BBMS}.
\newpage
\section*{Acknowledgments}

The author has been supported by MICIU/AEI/10.13039/501100011033 and ERDF/EU through the grant PID2021-122126NB-C31, and by ``Maria de Maeztu'' Excellence Unit IMAG, reference CEX2020-001105-M funded by MICIU/AEI/10.13039/501100011033.

Thanks are given to Audrey Fovelle, Gilles Godefroy, Vladimir Kadets, Ruben Medina, Abraham Rueda, and Dirk Werner for reading a first version of this manuscript and giving interesting remarks, and to Carlo De Bernardi for pointing out the result of Vesel\'{y} commented in Remark~\ref{remarkVesely}. We also thank the anonymous referee for the carefully reading of the manuscript and for some suggestions which have improved the final version of it.

\end{document}